\documentclass[11pt]{article}

\title{Association schemes related to universally optimal configurations,
Kerdock codes and extremal Euclidean line-sets}
\author{Kanat Abdukhalikov
\thanks{Supported by Japan Society for the Promotion of Science}\\
Institute of Mathematics, \\
Pushkin Str 125, Almaty 050010, Kazakhstan\\
abdukhalikov@math.kz\\
Eiichi Bannai\\
Graduate School of Mathematics, Kyushu University, \\
Hakozaki 6-10-1, Higashi-ku, Fukuoka 812-8581, Japan\\
bannai@math.kyushu-u.ac.jp\\
Sho Suda\\
Graduate School of Mathematics, Kyushu University, \\
Hakozaki 6-10-1, Higashi-ku, Fukuoka 812-8581, Japan}

\date{ }

\usepackage{amssymb}

\begin{document}

\maketitle

\newcommand{\Fq}{\mbox{$\mathbb{F}_q$}}         
\newcommand{\Fqq}{\mbox{$\mathbb{F}_q^*$}}      
\newcommand{\Zf}{\mbox{$\mathbb{Z}_4$}}         
\newcommand{\C}{\mbox{$\mathbb{C}$}}            

\newtheorem{theorem}{Theorem}
\newtheorem{lemma}[theorem]{Lemma}
\newtheorem{proposition}[theorem]{Proposition}

\begin{abstract}
H.~Cohn et. al. proposed an association scheme of 64 points in $\mathbb{R}^{14}$
which is conjectured to be a universally optimal code. We show that this scheme
has a generalization in terms of Kerdock codes, as well as in terms of maximal
real mutually unbiased bases. These schemes also related to extremal line-sets
in Euclidean spaces and Barnes-Wall lattices. D.~de Caen and E.~R.~van~Dam
constructed two infinite series of formally dual 3-class association schemes.
We explain this formal duality by constructing two dual abelian schemes related
to quaternary linear Kerdock and Preparata codes.

\end{abstract}


\section{Introduction}

In \cite{Cohn} Henry Cohn and others defined and studied universally optimal
configurations in Euclidean spaces, spherical point configurations
that minimize broad class of functions like potential energy functions.
They also obtained experimental results and conjectured
that  a some three class association scheme on 64 points determines universally
optimal configuration in $\mathbb{R}^{14}$.
This scheme is uniquely determined by their parameters \cite{Ban}
and has automorphism group $4^3:(2\times L_3(2))$, where $2\times L_3(2)$
is the stabilizer of a point.
It has the following first and second eigenmatrices:

$$
P=Q=\left(\begin{array}{rrrr}
         1   &  14  & 42  &  7  \\
         1   &  -6  &  6  & -1   \\
         1   &   2  & -2  & -1   \\
         1   &  -2  & -6  &  7
      \end{array}
\right) .
$$

The scheme generates a configuration of 64 vectors on a sphere of squared
radius 7 in $\mathbb{R}^{14}$. These vectors generate integral lattice
with automorphism group $2^{14}:(2^3:L_3(2))$.
The lattice can be obtained by construction A from binary
shortened projective [14,4,7] code.
Theta series of the lattice is equal to $\theta(q)=1+28q^4+1024q^7+2156q^8+\cdots$.
So the 64 vectors are not minimal vectors of the lattice.

We show that the scheme has a generalization in terms of binary and quaternary Kerdock
and Preparata codes, as well as in terms of maximal real mutually unbiased bases (MUB).
Starting from doubly shortened binary (shortened quaternary) Kerdock code of length $N$
(resp. of length $N/2$), where $N=2^{m+1}$ with odd $m$, one can construct a 3-class
association scheme of size $N^2/4$, which leads to a spherical code in
$\mathbb{R}^{N-2}$ of size $N^2/4$.
As binary Kerdock code we understand a binary code of length $N$ obtained
from a Kerdock set \cite{Kan}, so Kerdock code is inside of a second order
Reed-Muller code (see for details Section \ref{kerdock}), and it has nonzero
distances $2^m \pm 2^{(m-1)/2}$, $2^m$, $2^{m+1}$.
Furthermore, we define Kerdock-like code as a binary code of length $N$ with
nonzero distances $(N \pm \sqrt{N})/2$, $N/2$ and $N$
(in particular, Kerdock codes are Kerdock-like codes).
It is an open question whether Kerdock-like code is actually a Kerdock code.
We do not even know whether $N$ must have the form $2^{m+1}$.

Any collection of maximal mutually unbiased bases in $\mathbb{R}^N$
is in one-to-one correspondence with Kerdock-like codes.
In particular, any Kerdock code determines maximal MUB.
We show that in fact any maximal MUB in $\mathbb{R}^N$ (equivalently, any Kerdock-like
code) determines a 3-class association scheme in $\mathbb{R}^{N-2}$ of size $N^2/4$
with the same parameters as schemes obtained from Kerdock codes (Theorem \ref{theorem2}).

D. de Caen and E. R. van Dam \cite{Caen} constructed two infinite series of
formally dual 3-class association schemes, related to Kerdock sets.
We explain this formal duality by constructing two dual abelian schemes
related to quaternary linear Kerdock and Preparata codes (Theorem \ref{theorem1}).
The situation is similar to one, when the formal duality between binary nonlinear
Kerdock and Preparata codes was explained by duality between quaternary linear
Kerdock and Preparata codes \cite{Ham}.

We also note that maximal real mutually unbiased bases determine a 4-class
association scheme of size $N^2+2N$ in $\mathbb{R}^N$, as well as a 3-class
association scheme of size $N^2/2$ in $\mathbb{R}^{N-1}$ (which corresponds
to the association scheme obtained from a shortened Kerdock-like code).

R. L. Griess Jr. showed \cite{Gri} that the 64 point code and other tricosine
codes can also be constructed using minimal vectors of Barnes-Wall lattice.
We give later explanation of this phenomenon in terms of extremal line-sets and
Kerdock codes.

It seems that the schemes on shortened Kerdock codes also might be candidates for
being universally optimal (or optimal) configurations in $\mathbb{R}^{2^{m+1}-2}$.
We note that in $\mathbb{R}^{14}$ for $\theta$-code, $\theta=\frac{1}{7}$,
Levenshtein's bound is 69.6 (corresponding scheme has 64 points).
Similarly, in $\mathbb{R}^{62}$ for $\theta$-code, $\theta=\frac{3}{31}$,
Levenshtein's bound is 1081 (corresponding scheme has 1024 points).

The authors wish to thank H.~Cohn, R.~L.~Griess~Jr, O.~Musin, A.~Barg
and C.~Carlet for useful consultations.

\section{Quaternary Kerdock and Preparata codes and dual abelian
association schemes}

In this section we construct two dual abelian association schemes in terms of
quaternary linear Kerdock and Preparata codes.

By quaternary linear Kerdock and Preparata codes we mean
the class of codes determined in \cite{Cal}.
Recall that $\mathbb{Z}_4$-linear Kerdock code $K$ and
$\mathbb{Z}_4$-linear Preparata code $P$ are linear codes over $\mathbb{Z}_4$ of
length $q=2^m=N/2$, $m\geq 3$, $m$ odd:

$$ 0 \subset K \subseteq P \subset \mathbb{Z}_4^q .$$

They are dual codes: $K^{\perp}=P$. Moreover, $K=P$ for $m=3$ and $K\neq P$
for $m>3$. The image under the Gray map of the quaternary Kerdock (resp. Preparata)
code is binary nonlinear Kerdock (resp. "Preparata") code.
For $m=3$ the Gray image of $K=P$ is the famous binary nonlinear Nordstrom-Robinson
code of length 16. Consider shortened Kerdock and punctured Preparata codes:

$$ 0 \subset K_{\rm short} \subseteq P_{\rm punct} \subset \mathbb{Z}_4^{q-1} =A .$$

For a code $C$, the punctured code comes from deleting the
coordinate at position $i$, and the shortened code from
deleting the $0$ at $i$ from the words of the subcode of $C$
having 0's at $i$.
Since the automorphism group of  Kerdock code acts transitively
on coordinates we can consider shortening and puncturing at any
fixed (same) position.
Note that $K_{\rm short}^{\perp}= P_{\rm punct}$.
Therefore, one has nondegenerate bilinear pairing
$$(K_{\rm short}, A/P_{\rm punct}) \rightarrow \mathbb{Z}_4,$$
which gives us duality between $K_{\rm short}$ and $A/P_{\rm punct}$.
We have $|K_{\rm short}| = |A/P_{\rm punct}| = q^2 = 4^m$.
We can consider $A/P_{\rm punct}$ as character group of
$K_{\rm short}$ or vice versa. Therefore, an abelian association
scheme on $K_{\rm short}$ defines dual abelian scheme on
$A/P_{\rm punct}$ (cosets of $P_{\rm punct}$).

Shortened $\mathbb{Z}_4$-Kerdock code \cite{Cal,Ham} is a code
of length $2^m-1$, $m$ odd. It has $4^m=N^2/4$ codewords and nonzero
codewords have Lee weights $2^m+2^{(m-1)/2}$, $2^m-2^{(m-1)/2}$, and $2^m$.
We are going to show that the following relations on the shortened
Kerdock code will determine an abelian 3 class association scheme:

\begin{equation}
\label{kerd}
(x,y)\in\left\{\begin{array}{ll}
         R_0,   &  {\rm if} \ x-y \ {\rm has \ weight}\ 0,  \\
         R_1,   &  {\rm if} \ x-y \ {\rm has \ weight}\ 2^m+2^{(m-1)/2},  \\
         R_2,   &  {\rm if} \ x-y \ {\rm has \ weight}\ 2^m-2^{(m-1)/2},  \\
         R_3,   &  {\rm if} \ x-y \ {\rm has \ weight}\ 2^m.
        \end{array}
\right.
\end{equation}

Cosets of punctured $\mathbb{Z}_4$-Preparata code $C=P_{\rm punct}$
have Lee weights 0, 1 and 2. Furthermore, for
cosets $a + C$ we can choose $a=(0,\dots,0,\pm 1,0,\dots,0)$,
$a=(0,\dots,0, +1,\dots,-1,0,\dots,0)$ or
$a=(0,\dots,0,2,0,\dots,0)$ (see Lemma \ref{lem}).
The following relations

\begin{equation}
\label{prep}
(x,y)\in\left\{\begin{array}{ll}
         R'_0,   &  {\rm if} \ x-y = C,  \\
         R'_1,   &  {\rm if} \ x-y =(0,\dots,0,\pm 1,0,\dots,0) +C,  \\
         R'_2,   &  {\rm if} \ x-y =(0,\dots,0,+1,\dots,-1,0,\dots,0) +C,  \\
         R'_3,   &  {\rm if} \ x-y =(0,\dots,0,2,0,\dots,0) +C.  \\
        \end{array}
\right .
\end{equation}
on $A/C$ will define a three class association
scheme which is dual to the previous scheme.

\begin{theorem}
\label{theorem1}
The relations (\ref{kerd}) on codewords of shortened $\mathbb{Z}_4$-Kerdock code
define a three class (abelian) association scheme,
with the first and the second eigenmatrices  given by:

$$
P=\left(
  \begin{array}{rccr}
  1 & \frac{(N-2\sqrt{N})(N-2)}{8} & \frac{(N+2\sqrt{N})(N-2)}{8} & \frac{N}{2}-1 \\
  1 & -\frac{\sqrt{N}(N-4)}{8}     & \frac{\sqrt{N}(N-4)}{8}      & -1 \\
  1 & \frac{\sqrt{N}}{2}           & -\frac{\sqrt{N}}{2}          & -1 \\
  1 & -\frac{N-2\sqrt{N}}{4}       & -\frac{N+2\sqrt{N}}{4}       & \frac{N}{2}-1 \\
  \end{array}
  \right)
,
$$
$$
Q=\left(\begin{array}{rrrr}
 1 & N-2         & \frac{(N-2)(N-4)}{4} & \frac{N}{2}-1 \\
 1 & -\sqrt{N}-2 & \sqrt{N}+2           & -1 \\
 1 & \sqrt{N}-2  & -\sqrt{N}+2          & -1 \\
 1 &             -2 & -\frac{N}{2}+2    & \frac{N}{2}-1 \\
 \end{array}
 \right).
$$
The relations (\ref{prep}) on cosets of punctured $\Zf$-Preparata
code define an association scheme which is dual to the former scheme,
so the scheme has the following first and second eigenmatrices:
$P'=Q$, $Q'=P$.
\end{theorem}

If we take the particular class of Kerdock and Preparata codes
considered in \cite{Ham} then the scheme has
automorphism group $4^m:Aut(K_{\rm short})$, where
$Aut(K_{\rm short}) \cong 2\times L_3(2)$ for $m=3$ and
$Aut(K_{\rm short}) \cong 2\times (\mathbb{F}_{2^m}^*:Aut(\mathbb{F}_{2^m}))$
for $m>3$. We also note that the abelian group $K_{\rm short}$ is isomorphic to the
Galois ring $GR(4,m)=\mathbb{Z}_4[\xi]$, $\xi^{q-1}=1$, $q=2^m$.
Isomorphism is given by map
$\gamma \mapsto (Tr(\gamma\xi^0),Tr(\gamma\xi^1),\dots, Tr(\gamma\xi^{q-2}))$,
where $\gamma \in GR(4,m)$
(see for details \cite{Ham}).

\medskip

At first we study the structure of cosets of punctured $\mathbb{Z}_4$-Preparata
code $C=P_{\rm punct}$.

\begin{lemma}
\label{lem}
a) There exists a partition $A/C=V_0 \cup V_1 \cup V_2\cup V_3$
into four sets, where $V_0$, $V_1$, $V_2$ and $V_3$ are
cosets of the form $C$, $(0,\dots,0,\pm 1,0,\dots,0)+C$,
$(0,\dots,0, +1,0,\dots,0,-1,0,\dots,0)+C$
and $(0,\dots,2,\dots,0)+C$, respectively.

b) The previous statement remains true if the set $V_2$ is
enumerated by elements of the form
$(0,\dots,0, 1,0,\dots,0, 1,0,\dots,0)+C$
and $(0,\dots,0, -1,0,\dots,0, -1,0,\dots,0)+C$.

c) Finally, $V_2$ can be enumerated by elements of the form
$(0,\dots,0, 2,0,\dots,0, 1,0,\dots,0)+C$.
\end{lemma}

{\em Proof}. a) First we note that
$|V_0|=1$, $|V_1|=2(q-1)$, $|V_2|=(q-1)(q-2)$, $|V_3|=q-1$, so
$|V_0|+|V_1|+|V_2|+|V_3| =q^2=|A/C|$. It remains for us to prove
that elements of all $V_i$ are different. Minimum Lee weight of
$\mathbb{Z}_4$-Preparata
code $P$ is 6, therefore minimum Lee weight of punctured code $C$
is either 4 (if there is an element
$(2,\pm 1,\pm 1,\pm 1,\pm 1,0,\dots,0)\in P$) or 5.
If the elements of $V_i$ are not different, then $C$ contains
elements of the form $(1,1,-1,-1,0,\dots,0)$ or
$(2,1,-1,0,\dots,0)$, which means $(2,1,1,-1,-1,0,\dots,0)\in P$
or $(2,2,1,-1,0,\dots,0)\in P$.
These cases are not possible, since codewords of $P$ are zero
sum vectors ($(1,\dots,1)\in K$ by definition and $P=K^\perp$) and
the codeword $(0,0,2,2,0,\dots,0)$ of weight 4 does not belong to $P$.

Parts b) and c) can be proved similarly. 
$\Box$

\medskip

{\em Proof of Theorem \ref{theorem1}}. First we prove that the relations
(\ref{prep}) on cosets of punctured $\Zf$-Preparata code define an
association scheme. Consider the following elements
in $\C[A/C]$ (the group ring of the group $A/C$
over the field of complex numbers):

$$D_i = \sum_{v\in V_i}Z^v, \ 0\leq i\leq 3.$$

We will show that the subalgebra $<D_0,D_1,D_2,D_3>$, generated by
elements $D_0$, $D_1$, $D_2$ and $D_3$, is Schurian.

The following equalities are obtained just from counting elements in $D_i$:

$$D_1\cdot (D_0+D_1+D_2+D_3) = 2(q-1)(D_0+D_1+D_2+D_3),$$
$$D_2\cdot (D_0+D_1+D_2+D_3) = (q-1)(q-2)(D_0+D_1+D_2+D_3),$$
$$D_3\cdot (D_0+D_1+D_2+D_3) = (q-1)(D_0+D_1+D_2+D_3),$$

Further, we have
$$(D_0+D_3)^2 = (D_0+D_3)(D_0+D_3) =q(D_0+D_3),$$
which implies
$$D_3^2 = (q-1)D_0 + (q-2)D_3.$$

Lemma \ref{lem} b) implies that
$$D_1^2 = 2(q-1)D_0 + 4D_2 + 2D_3.$$

Finally, by lemma \ref{lem} c) we have
$$D_3D_1 = D_1 + 2D_2.$$

Matrices of multiplications by $D_1$, $D_2$, $D_3$ with respect to basis
$D_0$, $D_1$, $D_2$, $D_3$ are given by

$$
\rho_1=\left(\begin{array}{rrrr}
         0   &  2q-2  &     0    &  0  \\
         1   &     0  &  2(q-2)  &  1   \\
         0   &     4  &  2(q-4)  &  2   \\
         0   &     2  &  2(q-2)  &  0
      \end{array}
\right) ,
$$
$$
\rho_2=\left(\begin{array}{rrrr}
         0   &      0   &  (q-1)(q-2)  & 0  \\
         0   &  2(q-2)  &  (q-4)(q-2)  & q-2   \\
         1   &  2(q-4)  &   q^2-6q+12  & q-3   \\
         0   &  2(q-2)  &  (q-3)(q-2)  & 0
               \end{array}
\right) ,
$$
$$
\rho_3=\left(\begin{array}{rrrr}
         0   &  0  &    0  &  q-1  \\
         0   &  1  &  q-2  &  0   \\
         0   &  2  &  q-3  &  0   \\
         1   &  0  &    0  &  q-2
      \end{array}
\right) .
$$

It is easy to see that  vectors $v_1=(1, 2q-2, (q-1)(q-2), q-1)$,
$v_2=(1, -\sqrt{2q}-2, \sqrt{2q}+2, -1)$,
$v_3=(1,\sqrt{2q}-2, -\sqrt{2q}+2, -1)$, $v_4=(1, -2, -q+2, q-1)$
are common left eigenvectors for matrices $\rho_i={^tB_i}$.
They are strokes of the matrix
$$
P'=Q=\left(\begin{array}{rrrr}
         1   &  2q-2         & (q-1)(q-2)    & q-1  \\
         1   & -\sqrt{2q}-2  & \sqrt{2q}+2   & -1   \\
         1   &  \sqrt{2q}-2  & -\sqrt{2q}+2  & -1   \\
         1   &  -2           & -q+2          & q-1
      \end{array}
\right).
$$

Now we are going to prove that relations (\ref{kerd}) define
an association scheme dual to scheme (\ref{prep}).
The latter abelian scheme determines a dual scheme on
$K_{\rm short}$, with partition
$K_{\rm short} = V'_0 \cup V'_1 \cup V'_2\cup V'_3$.
We prove that this partition corresponds to sets of codewords
of Lee weights 0, $2^m+2^{(m-1)/2}$, $2^m-2^{(m-1)/2}$, and $2^m$
respectively.

Indeed, according to \cite[section 4.7.1]{Cam}, for any element
$u\in V'_j$ we should have

$$Q_{jk} = \sum_{v\in V_k}i^{-(v,u)},$$
where $i=\sqrt{-1}$. Such equations determine elements in $V'_j$
uniquely. For example, let us take codeword $u\in K_{\rm short}$
of Lee weight $2^m+2^{(m-1)/2}$ (other cases can be considered
analogously). Then codeword $u$ has form
$(2^a, 1^b,(-1)^c, 0^{q-1-a-b-c})$, $2a+b+c=2^m + 2^{(m-1)/2}$.
Since $2u=(0^a,2^b,2^c,0^{q-1-a-b-c})\in K_{\rm short}$,
we have $b+c=2^{m-1}$ and $a=(q+\sqrt{2q})/4$.

For $k=1$ we have $V_1=\{(0,\dots,0,\pm 1,0,\dots,0)+C\}$,
$|V_1|=2(q-1)$, and
\begin{eqnarray*}
\sum_{v\in V_1}i^{-(v,u)} & = &
         \frac{q+\sqrt{2q}}{4}i^2 + bi^{-1} + ci + (\frac{q-\sqrt{2q}}{4}-1) + \\
                          &   &
         \frac{q+\sqrt{2q}}{4}i^2 + bi + ci^{-1} + (\frac{q-\sqrt{2q}}{4}-1) = -\sqrt{2q}-2.
\end{eqnarray*}
Therefore, $u\in V'_1$.

Finally we show that one can use a shortening (puncturing) at any position.
It is enough to show that the automorphism group of Kerdock code is transitive
on coordinates. It just follows from the definition of $\mathbb{Z}_4$-Kerdock
code \cite{Cal}. Let $V$ be a $m$-dimensional vector space over $\mathbb{F}_2$ and
$R$ be a binary symmetric $m\times m$ matrix. The map $T_R: V\rightarrow \mathbb{Z}_4$
is given by
$$T_R(v)=\sum_{i=1}^m R_{ii}\widehat{v}_i^2 + 2\sum_{i<j}R_{ij}\widehat{v}_i\widehat{v}_j ,$$
where the entries of the matrix $R$ are identified with the elements 0, 1
of $\mathbb{Z}_4$ and $\widehat{v}\in \mathbb{Z}_4^m$ is congruent to $v\in V$
modulo 2
(such map is called $\mathbb{Z}_4$-valued quadratic form on $V$).
It is easy to see that

\begin{equation}
\label{identity}
T_R(u+v)=T_R(u)+T_R(v)+2\widehat{u}R\widehat{v}^T.
\end{equation}

Let $S$ be a set of $2^m$ binary symmetric $m\times m$ matrices
such that the difference of any two is non-singular. Then the
associated $\mathbb{Z}_4$-Kerdock code is defined as the set
$$\{(T_R(v)+2\widehat{s}\cdot \widehat{v} +\varepsilon)_v \mid
R\in S, \ s\in V, \ \varepsilon\in \mathbb{Z}_4\}$$
Now (\ref{identity}) implies that the automorphism group of a
Kerdock code is transitive on coordinates, since each translation
$v \mapsto v+u$ of $V$ is an automorphism of Kerdock code.

Theorem is proved.
$\Box$

\section{Binary Kerdock codes and association schemes}
\label{kerdock}

In this section we discuss binary Kerdock codes and related association schemes.

Recall the definition of binary nonlinear Kerdock codes \cite{Cal}.
Let $V$ be a $(m+1)$-dimensional vector space over $\mathbb{F}_2$, $m$ odd.
Any polynomial function $f$ on $V$ can be considered as codeword of length
$2^{m+1}$, evaluating $f$ on all vectors of $V$:
$(f(v_1),\dots,f(v_{2^{m+1}}))$.
The set of all polynomials of degree at most $k$ is called $k$th order
Reed-Muller code $RM(k,m+1)$. The first order Reed-Muller code
$RM(1,m+1)$ consists of affine functions.

The second order Reed-Muller code $RM(2,m+1)$ consists of elements $Q+f$,
where $f\in RM(1,m+1)$ and $Q$ is a quadratic form.
A quadratic form is a map $Q: V \mapsto \mathbb{F}_2$ such that
$Q(0)=0$ and
$$B_Q(x,y)=Q(x+y)+Q(x)+Q(y)$$
is bilinear.
For quadratic form $Q$ and affine function $f$ the element $Q+f \in RM(2,m+1)$
has $2^m \pm 2^{(m-1+r)/2}$ zeroes, where $r$ is the dimension of the radical
of the bilinear form $B_Q$, associated with $Q$.

A quadratic form is called non-singular if its associated bilinear form
is non-singular. Non-singular form $Q$ has a type $\chi(Q)=\pm 1$, where
$Q$ has precisely $2^m +\chi(Q)2^{(m-1)/2}$ zeroes.
The projective quadrics associated to non-singular quadratic forms with
$\chi(Q)= +1$ are called hyperbolic, those with $\chi(Q)=-1$ elliptic.

A Kerdock set is a set of $2^m$ alternating bilinear forms
(or, equivalently,  binary skew-symmetric $(m+1)\times(m+1)$ matrices)
such that the difference of any two is non-singular.
Let $S$ be a Kerdock set. Then the binary Kerdock code $\mathcal{K}(S)$ is
a set of elements $Q+f$, where $f\in RM(1,m+1)$ and
$Q$ is a quadratic form such that its associated bilinear form $B_Q$ belongs to $S$.
So Kerdock code $\mathcal{K}$ lives in $RM(2,m+1)$ and contains $RM(1,m+1)$:
$$RM(1,m+1) \subset \mathcal{K} \subset RM(2,m+1).$$

Actually the quadratic forms, generating a Kerdock code, are chosen so
that the minimum distance between any of two of the cosets $Q+RM(1,m+1)$
and $Q'+RM(1,m+1)$ is as large as possible, namely $2^m - 2^{(m-1)/2}$,
which occurs if and only if $Q+Q'$ is non-singular;
and Kerdock set provides the maximum possible number of such cosets.

\medskip

In \cite{Cal} it was shown that any binary Kerdock code can be obtained
from a $\mathbb{Z}_4$-Kerdock code with the aid of Gray map.
If this code is $\mathbb{Z}_4$-linear then its shortened code
defines an association scheme, as we showed in the previous section.
Actually the structure of an association scheme can be defined for
any doubly shortened binary Kerdock code and, therefore, for any shortened
$\mathbb{Z}_4$-Kerdock code.

\begin{proposition}
\label{dshort}
Let $C$ be a shortened in two (arbitrary) positions binary Kerdock code.
Then distinct codewords in $C$ have only distances
$2^m+2^{(m-1)/2}$, $2^m-2^{(m-1)/2}$, $2^m$
and these values determine three class association scheme on $q^2=N^2/4$
points with the first and second eigenmatrices $P$ and $Q$.
\end{proposition}

{\em Proof}. It follows from Theorem 8 in \cite{Caen}.
There it was constructed five class association schemes, and by
fusion ($R_1\cup R_2$, $R_3\cup R_4$, $R_5$) one can get a three
class schemes with required parameters.
We need to show that vertices of that scheme are actually
codewords of doubly shortened Kerdock code. Indeed, the vertex set consists
of all ordered pairs $(B,Q)$, where $Q$ is a quadratic form,
such that its associated bilinear form $B$ belongs to a fixed
Kerdock set $S$, and $Q(v)=0$ for a fixed vector in $V$.
Actually we can remove $B$ from the description of the vertex set,
since it can be determined uniquely from $Q$. Therefore
the vertex set consists of quadratic forms $Q$, which are
codewords from the Kerdock code, defined from Kerdock set $S$.
Moreover, we have $Q(0)=Q(v)=0$, so these codewords determine
doubly shortened (in $0$ and $v$) code.

Further, we note that if $Q\not\equiv Q' \pmod{RM(1,m+1)}$ then distance
$$d(Q,Q')= \# {\rm nonzeros} (Q+Q')=
2^{m+1}-(2^m +\chi(Q+Q')2^{(m-1)/2})=2^m -\chi(Q+Q')2^{(m-1)/2},$$
so
$$d(Q,Q')= 2^m \pm 2^{(m-1)/2} \Leftrightarrow \chi(Q+Q')=\mp 1.$$
If $Q\equiv Q' \pmod{RM(1,m+1)}$ then $d(Q,Q')= 2^m$.

Finally, we note that the automorphism group of a Kerdock code is
transitive on coordinates, since each translation
$u \mapsto u+a$ of $V$ is an automorphism of Kerdock code
(it leaves invariant each coset $Q+RM(1,m+1)$).
Using this automorphism and taking $v=b-a$ we can send the
pair $(0,v)$ to any fixed pair $(a,b)$.
$\Box$

\medskip

One might ask what happens if we consider a shortened binary Kerdock code.

\begin{proposition}
\label{short-dels}
Let $C$ be a shortened binary Kerdock code. Then distinct codewords
in $C$ have only distances $2^m+2^{(m-1)/2}$, $2^m-2^{(m-1)/2}$, $2^m$
and these values determine three class association scheme on $2q^2=N^2/2$ points
with the following first and second eigenmatrices:

$$
P=\left(\begin{array}{rccr}
  1 & \frac{(N-\sqrt{N})(N-2)}{4} & \frac{(N+\sqrt{N})(N-2)}{4} & N-1 \\
  1 & -\frac{\sqrt{N}(N-2)}{4} & \frac{\sqrt{N}(N-2)}{4} & -1 \\
  1 & \frac{\sqrt{N}}{2} & -\frac{\sqrt{N}}{2} & -1 \\
  1 & -\frac{N-\sqrt{N}}{2} & -\frac{N+\sqrt{N}}{2} & N-1 \\
  \end{array}
\right),
$$

$$
Q=\left(\begin{array}{rrrr}
 1 & N-1 & \frac{(N-2)(N-1)}{2} & \frac{N}{2}-1 \\
 1 & -\sqrt{N}-1 & \sqrt{N}+1 & -1 \\
 1 & \sqrt{N}-1 & -\sqrt{N}+1 & -1 \\
 1 & -1 & -\frac{N}{2}+1 & \frac{N}{2}-1 \\
 \end{array}
\right) .
$$
\end{proposition}

{\em Proof}. See Delsarte's thesis \cite{Del}, Example 2 on page 82.
$\Box$

\medskip

In his thesis \cite[Example 2 on page 82]{Del} Delsarte considers
a binary code $Y$ of length $N-1$ with three nonzero distances
$2^m \pm 2^{(m-1)/2}$, $2^m$. He shows that $|Y|\leq N^2/2$,
so shortened Kerdock code is optimal in the sense that it is a code
with maximal possible number of codewords
among binary codes of length  $N-1$ with mentioned three nonzero distances.
We repeated his reasoning for a doubly shortened Kerdock code.
Let us take a binary code  $Y$ of length $N-2$ with three nonzero distances
$2^m \pm 2^{(m-1)/2}$, $2^m$. Calculations show that
$$\alpha(z) =
   4N^{-2}|Y| (\Phi_0(z) + \frac{N+2}{2N-2} \Phi_1(z) +
   \frac{3}{N-1} \Phi_2(z) + \frac{3}{2N-2} \Phi_3(z)),$$
where notations are taken from \cite{Del}. Then theorem 5.23(ii) from \cite{Del} implies
that $|Y| \leq  N^2/4$.
It means that doubly shortened Kerdock code is a code
with maximal possible number of codewords
among binary codes of length  $N-2$ with three nonzero distances
$2^m \pm 2^{(m-1)/2}$, $2^m$.

\section{Mutually unbiased bases and association schemes}
\label{mut}

In the previous sections, in particular in Theorem \ref{theorem1}
and Proposition \ref{dshort},
we have constructed 3-class association schemes $Y$ in $\mathbb{R}^{N-2}$ of size
$|Y|=N^2/4$,  with $N=2^{m+1}$ and $m$ odd, from Kerdock codes.
These association schemes are generalizations of the association scheme
of 64 points in $\mathbb{R}^{14}$ which is a candidate of universal optimal code
considered by Cohn and Kumar \cite{Cohn2}. In this section we will show that from such Kerdock
codes we can obtain the specific line systems in real Euclidean space $\mathbb{R}^N$,
or equivalently, maximal real MUB (mutually unbiased bases) in real Euclidean space
$\mathbb{R}^N$. (The reader is referred to \cite{Boy} for the details on MUB.)
It is an interesting question whether any such line system in real Euclidean
space $\mathbb{R}^N$ can be obtained from a Kerdock code, but
this question is still open.
We will show also that from any maximal real MUB, or equivalently from any such
line system in real Euclidean space $\mathbb{R}^N$, we obtain an association scheme
which has the same parameters as those in the association scheme
obtained in Theorem \ref{theorem1}. So, this result (Theorem \ref{theorem2} below) is a
generalization of Proposition \ref{dshort}. Our proof is based on the
method due to Delsarte-Goethals-Seidel \cite{Del3}, and does not depend
on the specific structures of the codes.

Two orthonormal bases $\mathcal{B}$ and $\mathcal{B}'$ in $\mathbb{R}^N$
are called mutually unbiased if $|(x,y)|=1/\sqrt{N}$ for any
$x\in \mathcal{B}$ and $y\in \mathcal{B}'$.
It is known \cite{Cal,Del2} that there can be at most $N/2 +1$ mutually unbiased
bases in dimension $N$, and constructions reaching this upper bound are
known only for values $N=2^{m+1}$.
Let us define Kerdock-like codes as binary codes of length $N$ with
$N^2$ codewords and with
nonzero distances $(N \pm \sqrt{N})/2$, $N/2$ and $N$
(in particular, Kerdock codes are Kerdock-like codes).
Any collection of $N/2 +1$ mutually unbiased bases is in one-to-one
correspondence with Kerdock-like codes. Indeed, let $\mathcal{B}_0$,
$\mathcal{B}_1$, \dots, $\mathcal{B}_{N/2}$ be mutually unbiased bases,
$\mathcal{B}_0$ be the standard orthonormal basis.
Consider the set $L$ of all vectors $x \in \mathbb{R}^N$ such that $x$ or $-x$
belongs to one of bases $\mathcal{B}_1$, \dots, $\mathcal{B}_{N/2}$.
Then $|L|=N^2$. The elements of $L$ have form
$(\pm 1/\sqrt{N}, \dots, \pm 1/\sqrt{N})$, since $\mathcal{B}_0$
and $\mathcal{B}_i$ ($i>0$) are mutually unbiased. One can convert $L$ to a
binary code $C$ by changing $1/\sqrt{N}$ and $-1/\sqrt{N}$ to
$0$ and $1$, respectively. Then it is easy to see that $C$ is a Kerdock-like code
(the conditions $(x,y)=\pm 1/\sqrt{N}$, $0$, $1$ for $x$, $y\in L$ mean
exactly that distances $d(x',y')= (N\mp \sqrt{N})/2$, $N/2$, $N$ for corresponding
images $x'$, $y'\in C)$.

In the construction of association schemes we used doubly shortened
binary Kerdock codes. Now we do similar procedure for mutually unbiased bases:
take vectors from $L$ with $1/\sqrt{N}$ in two fixed coordinates and then
delete these coordinates.
We will get a configuration of $N^2/4$ vectors in $\mathbb{R}^{N-2}$ of equal
length, such that cosines of angles  between distinct vectors are equal to
$\frac{-\sqrt{N}-2}{N-2}$, $\frac{-2}{N-2}$ and $\frac{\sqrt{N}-2}{N-2}$.

\begin{theorem}
\label{theorem2}
Let $M$ be a maximal mutually unbiased bases in $\mathbb{R}^N$ and $X=M\cup (-M)$.
For $u,v\in X$ such that $(u,v)=0$, we put
$Y':=\{x\in X|(u,x)=(v,x)=\frac{1}{\sqrt{N}}\}$, and let
$Y$ be the orthogonal projection of vectors $Y'$ to $\langle u, v\rangle^\perp$
rescaled to make a spherical code.
Then $Y$ is an association scheme with the first and second eigenmatrices given by:
$$
P=\left(
  \begin{array}{rccr}
  1 & \frac{(N-2\sqrt{N})(N-2)}{8} & \frac{(N+2\sqrt{N})(N-2)}{8} & \frac{N}{2}-1 \\
  1 & -\frac{\sqrt{N}(N-4)}{8}     & \frac{\sqrt{N}(N-4)}{8}      & -1 \\
  1 & \frac{\sqrt{N}}{2}           & -\frac{\sqrt{N}}{2}          & -1 \\
  1 & -\frac{N-2\sqrt{N}}{4}       & -\frac{N+2\sqrt{N}}{4}       & \frac{N}{2}-1 \\
  \end{array}
  \right),
  $$

$$
Q=\left(\begin{array}{rrrr}
 1 & N-2         & \frac{(N-2)(N-4)}{4} & \frac{N}{2}-1 \\
 1 & -\sqrt{N}-2 & \sqrt{N}+2           & -1 \\
 1 & \sqrt{N}-2  & -\sqrt{N}+2          & -1 \\
 1 & -2          & -\frac{N}{2}+2       & \frac{N}{2}-1 \\
 \end{array}
 \right).
$$
\end{theorem}

{\em Proof}.
We put $\alpha=\frac{-\sqrt{N}-2}{N-2}$, $\beta=\frac{\sqrt{N}-2}{N-2}$,
$\gamma=\frac{-2}{N-2}$.
Since $A(Y)=\{\alpha, \beta, \gamma \}$, Y is a 3-distance set.
The annihilator polynomial
$F(x):=\prod_{\alpha \in A(Y)}\frac{x-\alpha}{1-\alpha}$
has the Gegenbauer polynomial expansion
$$F(x)=\frac{4}{N^2}Q_0(x)+\frac{2(N^2+6)(N-2)}{N^3(N-1)}Q_1(x)+
\frac{(N-2)^3(N+3)}{N^3(N-1)}Q_2(x)+\frac{6(N-2)(N-3)}{N^2(N-1)}Q_3(x).$$

As $|Y|=N^2/4$ , Theorem 6.5 of \cite{Del3} implies that Y is a spherical 3-design.

By Lemma 7.3 of \cite{Del3}, for $0 \leq i, j \leq 2$, $i+j \neq 4$ and $z :=(\xi,\eta)$,
the intersection numbers $p_{\alpha,\beta}(\xi,\eta)$ satisfy the linear equation
$$ \sum_{x ,y \in A(Y)}x^i y^j p_{x,y}(\xi,\eta)=
\frac{N^2}{4}F_{i,j}(z )-z^i-z^j+\delta_{1,z},$$
where ${\displaystyle t^i=\sum_{k=0}^i f_{i,k} Q_k(t)}$
 and ${\displaystyle F_{i,j}(t)=\sum_{k=0}^{min(i,j)} f_{i,k} f_{j,k} Q_k(t)}$.

Now let $z =(\xi,\eta)$ be fixed.  Then we get following equation:

\begin{eqnarray*}
\left(
  \begin{array}{cccccccc}
       1&    1&    1&    1&    1&    1&    1&    1    \\
       \alpha &    \beta &   \gamma &    \alpha &    \beta &    \gamma &    \alpha &    \beta   \\
      \alpha^2 &  \beta^2  &  \gamma^2  &  \alpha^2  & \beta^2   & \gamma^2   &  \alpha^2  & \beta^2    \\
      \alpha &  \alpha  & \alpha   & \beta   & \beta   & \beta   & \gamma   & \gamma      \\
 \alpha^2 &  \alpha \beta  & \alpha \gamma  & \alpha \beta   & \beta^2   & \beta \gamma   &  \alpha  \gamma  & \beta \gamma   \\
 \alpha^3 &  \alpha \beta^2  & \alpha \gamma^2  & \alpha^2 \beta   & \beta^3   & \beta \gamma^2   &  \alpha^2  \gamma  & \beta^2 \gamma     \\
 \alpha^2 &    \alpha^2 &   \alpha^2 &    \beta^2 &    \beta^2 &    \beta^2 &    \gamma^2 &    \gamma^2  \\
 \alpha^3 &    \alpha^2 \beta &   \alpha^2 \gamma &    \alpha \beta^2 &    \beta^3 &    \beta^2 \gamma &   \alpha \gamma^2 &  \beta  \gamma^2   \\
  \end{array}
\right)
\left(
  \begin{array}{c}
   p_{\alpha,\alpha}(\xi,\eta)    \\
   p_{\beta,\alpha}(\xi,\eta)    \\
   p_{\gamma,\alpha}(\xi,\eta)    \\
   p_{\alpha,\beta}(\xi,\eta)    \\
   p_{\beta,\beta}(\xi,\eta)    \\
   p_{\gamma,\beta}(\xi,\eta)    \\
   p_{\alpha,\gamma}(\xi,\eta)    \\
   p_{\beta,\gamma}(\xi,\eta)    \\
  \end{array}
\right) = \\
\left(
  \begin{array}{c}
   F_{0,0}(z )-z^0-z^0+\delta_{1,z} -\gamma^0  p_{\gamma,\gamma}(\xi,\eta)   \\
   F_{1,0}(z )-z^1-z^0+\delta_{1,z} -\gamma^1  p_{\gamma,\gamma}(\xi,\eta)   \\
   F_{2,0}(z )-z^2-z^0+\delta_{1,z} -\gamma^2  p_{\gamma,\gamma}(\xi,\eta)   \\
   F_{0,1}(z )-z^0-z^1+\delta_{1,z} -\gamma^1  p_{\gamma,\gamma}(\xi,\eta)   \\
   F_{1,1}(z )-z^1-z^1+\delta_{1,z} -\gamma^2  p_{\gamma,\gamma}(\xi,\eta)   \\
   F_{2,1}(z )-z^2-z^1+\delta_{1,z} -\gamma^3  p_{\gamma,\gamma}(\xi,\eta)   \\
   F_{0,2}(z )-z^0-z^2+\delta_{1,z} -\gamma^2  p_{\gamma,\gamma}(\xi,\eta)   \\
   F_{1,2}(z )-z^1-z^2+\delta_{1,z} -\gamma^3  p_{\gamma,\gamma}(\xi,\eta)   \\
  \end{array}
\right),
\end{eqnarray*}
where
$$
\det
\left(
  \begin{array}{cccccccc}
       1&    1&    1&    1&    1&    1&    1&    1    \\
       \alpha &    \beta &   \gamma &    \alpha &    \beta &    \gamma &    \alpha &    \beta   \\
      \alpha^2 &  \beta^2  &  \gamma^2  &  \alpha^2  & \beta^2   & \gamma^2   &  \alpha^2  & \beta^2    \\
      \alpha &  \alpha  & \alpha   & \beta   & \beta   & \beta   & \gamma   & \gamma      \\
 \alpha^2 &  \alpha \beta  & \alpha \gamma  & \alpha \beta   & \beta^2   & \beta \gamma   &  \alpha  \gamma  & \beta \gamma   \\
 \alpha^3 &  \alpha \beta^2  & \alpha \gamma^2  & \alpha^2 \beta   & \beta^3   & \beta \gamma^2   &  \alpha^2  \gamma  & \beta^2 \gamma     \\
 \alpha^2 &    \alpha^2 &   \alpha^2 &    \beta^2 &    \beta^2 &    \beta^2 &    \gamma^2 &    \gamma^2  \\
 \alpha^3 &    \alpha^2 \beta &   \alpha^2 \gamma &    \alpha \beta^2 &    \beta^3 &    \beta^2 \gamma &   \alpha \gamma^2 &  \beta  \gamma^2   \\
  \end{array}
\right)
=(\alpha-\beta)^6(\alpha-\gamma)^4(\beta-\gamma)^4 \neq 0.
$$

In case of $(\xi,\eta)=1 $ we have $p_{\gamma,\gamma}(\xi,\eta)= N/2 -1$.
In case of $(\xi,\eta)=\alpha $ or $ \beta$, we have $p_{\gamma,\gamma}(\xi,\eta)=0$.
Finally, in case of $(\xi,\eta)=\gamma$ we have $p_{\gamma,\gamma}(\xi,\eta)= N/2 -2$.

Therefore,  intersection numbers $p_{x,y}(\xi,\eta)$ are determined uniquely by
$p_{\gamma,\gamma}(\xi,\eta)$. Hence Y is an association scheme.
This completes the proof of Theorem \ref{theorem2}. $\Box$

\medskip
{\bf Remark.} In Theorem 7.4 of  \cite{Del3} it is mentioned that if
$t\geq 2s-3$, then for any fixed $z=\tilde{\gamma} =(\xi,\eta)$,
the intersection numbers $p_{\tilde{\alpha},\tilde{\beta}}(\xi,\eta)$
are uniquely determined by $p_{\tilde{\gamma},\tilde{\gamma}}(\xi,\eta)$.
Our claim is that intersection numbers
$p_{\tilde{\alpha},\tilde{\beta}}(\xi,\eta)$ are uniquely determined
by $p_{\gamma,\gamma}(\xi,\eta)$ with a suitable $\gamma$
which is not necessarily equal to $\tilde{\gamma}$.  So, our argument
is slightly general than in Theorem 7.4 of  \cite{Del3}.

\medskip
Here we mention some more useful information. We can easily
calculate the intersection matrices
$B_i=(p_{i,j}^k)_{0\leq j \leq 3 , 0\leq k \leq 3 }$
and Krein parameter matrices
$B_i^*=(q_{i,j}^k)_{0\leq j \leq 3 , 0\leq k \leq 3 }$.
Namely, they are given as follows:

$B_1=\left (\begin{array}{cccc}
 0 & 1 & 0 & 0  \\
\frac{(N-2\sqrt{N})(N-2)}{8}  & \frac{(N+2\sqrt{N})(N-7\sqrt{N}+12)}{16} & \frac{(N-2\sqrt{N})(N-\sqrt{N}-4)}{16} & \frac{(N-2\sqrt{N})(N-2\sqrt{N}-4)}{16}\\
0 & \frac{(N+2\sqrt{N})(N-\sqrt{N}-4)}{16} & \frac{(N-2\sqrt{N})(N+\sqrt{N}-4)}{16} & \frac{N(N-4)}{16}   \\
0 & \frac{N-2\sqrt{N}-4}{4} & \frac{N-2\sqrt{N}}{4} & 0 \\
\end{array}
\right)$, \\

$B_2=\left (\begin{array}{cccc}
0 & 0 & 1 & 0 \\
0 & \frac{(N+2\sqrt{N})(N-\sqrt{N}-4)}{16} & \frac{(N-2\sqrt{N})(N+\sqrt{N}-4)}{16} & \frac{N(N-4)}{16}   \\
\frac{(N+2\sqrt{N})(N-2)}{8} & \frac{(N+2\sqrt{N})(N+\sqrt{N}-4)}{16} & \frac{(N-2\sqrt{N})(N+7\sqrt{N}+12)}{16} & \frac{(N+2\sqrt{N})(N+2\sqrt{N}-4)}{16} \\
0 & \frac{N+2\sqrt{N}}{4} & \frac{N+2\sqrt{N}-4}{4} & 0 \\
\end{array}
\right)$, \\

$B_3=\left (\begin{array}{cccc}
0 & 0 & 0 & 1 \\
0 & \frac{N-2\sqrt{N}-4}{4} & \frac{N-2\sqrt{N}}{4} & 0 \\
0 & \frac{N+2\sqrt{N}}{4} & \frac{N+2\sqrt{N}-4}{4} & 0 \\
\frac{N}{2}-1 & 0 & 0 & \frac{N}{2}-2 \\
\end{array}
\right)$ ;\\

$B_1^*=\left (\begin{array}{cccc}
 0 & 1 & 0 & 0  \\
 N-2 & 0 & 4 & 2 \\
 0 & N-4 & N-8 & N-4  \\
 0 & 1 & 2 & 0 \\
\end{array}
\right)$ ,\\

$B_2^*=\left (\begin{array}{cccc}
 0 & 0 & 1 & 0  \\
 0 & N-4 & N-8 & N-4 \\
 \frac{(N-2)(N-4)}{4} & \frac{(N-4)(N-8)}{4} & \frac{N^2-12N+48}{4} & \frac{(N-4)(N-6)}{4}  \\
 0 & \frac{N}{2}-2 & \frac{N}{2}-3 & 0 \\
\end{array}
\right)$, \\

$B_3^*=\left (\begin{array}{cccc}
 0 & 0 & 0 & 1  \\
 0 & 1 & 2 & 0 \\
 0 & \frac{N}{2}-2 & \frac{N}{2}-3 & 0  \\
 \frac{N}{2}-1 & 0 & 0 & \frac{N}{2}-2 \\
\end{array}
\right)$. \\

\bigskip

The original association scheme $X$ of size $N^2+2N$ in $\mathbb{R}^N$
(attached to a real MUB) is a class 4 association scheme with the following
parameters. Note that this is a $5$-spherical design and of degree $s=4$.
It is interesting to note that they are $Q$-polynomial association schemes
(and not $P$-polynomial association scheme for $N\geq 4 $).
The reader is referred to Bannai-Bannai \cite{Ban2} for more details,
where this fact was first noticed. These association schemes
are possible candidates of universally optimal codes in the sense
of Cohn-Kumar \cite{Cohn2}.

$B_1=\left (\begin{array}{ccccc}
 0 & 1 & 0 & 0 & 0 \\
\frac{N^2}{2}  & \frac{(N+\sqrt{N})(N-2)}{4} & \frac{N^2}{4} & \frac{(N-\sqrt{N})(N-2)}{4} & 0 \\
0 & N-1 & 0 & N-1 & 0   \\
0 & \frac{(N-\sqrt{N})(N-2)}{4} & \frac{N^2}{4} & \frac{(N+\sqrt{N})(N-2)}{4} & \frac{N^2}{2} \\
0 & 0 & 0 & 1 & 0 \\
\end{array}
\right)$ ,\\

$B_2=\left (\begin{array}{ccccc}
0 & 0 & 1 & 0 & 0 \\
0 & N-1 & 0 & N-1 & 0 \\
2(N-1) & 0 & 2(N-2) & 0 & 2(N-1) \\
0 & N-1 & 0 & N-1 & 0 \\
0 & 0 & 1 & 0 & 0 \\
\end{array}
\right)$ ,\\

$B_3=\left (\begin{array}{ccccc}
0 & 0 & 0 & 1 & 0 \\
0 & \frac{(N-\sqrt{N})(N-2)}{4} & \frac{N^2}{4} & \frac{(N+\sqrt{N})(N-2)}{4} & \frac{N^2}{2} \\
0 & N-1 & 0 & N-1 & 0 \\
\frac{N^2}{2} & \frac{(N+\sqrt{N})(N-2)}{4} & \frac{N^2}{4} & \frac{(N-\sqrt{N})(N-2)}{4} & 0 \\
0 & 1 & 0 & 0 & 0 \\
\end{array}
\right)$ ,\\

$B_4=\left (\begin{array}{ccccc}
0 & 0 & 0 & 0 & 1 \\
0 & 0 & 0 & 1 & 0 \\
0 & 0 & 1 & 0 & 0 \\
0 & 1 & 0 & 0 & 0 \\
1 & 0 & 0 & 0 & 0 \\
\end{array}
\right)$ ;\\

$B_1^*=\left (\begin{array}{ccccc}
 0 & 1 & 0 & 0 & 0 \\
 N & 0 & \frac{2N}{N+2} & 0 & 0 \\
 0 & N-1 & 0 & N-1 & 0  \\
 0 & 0 & \frac{N^2}{N+2} & 0 & N \\
 0 & 0 & 0 & 1 & 0 \\
\end{array}
\right)$ ,\\

$B_2^*=\left (\begin{array}{ccccc}
 0 & 0 & 1 & 0 & 0 \\
 0 & N-1 & 0 & N-1 & 0 \\
 \frac{(N+2)(N-1)}{2} & 0 & \frac{(N+2)(N-2)}{2} & 0 & \frac{(N+2)(N-1)}{2} \\
 0 & \frac{N(N-1)}{2} & 0 & \frac{N(N-1)}{2} & 0 \\
 0 & 0 & \frac{N}{2} & 0 & 0 \\
\end{array}
\right)$ ,\\

$B_3^*=\left (\begin{array}{ccccc}
 0 & 0 & 0 & 1 & 0 \\
 0 & 0 & \frac{N^2}{N+2} & 0 & N \\
 0 & \frac{N(N-1)}{2} & 0 & \frac{N(N-1)}{2} & 0  \\
 \frac{N^2}{2} & 0 & \frac{N^3}{2(N+2)} & 0 & \frac{N(N-2)}{2} \\
 0 & \frac{N}{2} & 0 & \frac{N}{2}-1 & 0 \\
\end{array}
\right)$ ,\\

$B_4^*=\left (\begin{array}{ccccc}
0 & 0 & 0 & 0 & 1 \\
0 & 0 & 0 & 1 & 0 \\
0 & 0 & \frac{N}{2} & 0 & 0 \\
0 & \frac{N}{2} & 0 & \frac{N}{2}-1 & 0 \\
\frac{N}{2} & 0 & 0 & 0 & \frac{N}{2}-1 \\
\end{array}
\right)$ ;\\

$P=\left(\begin{array}{ccccc}
  1 & \frac{N^2}{2} & 2(N-1) & \frac{N^2}{2} & 1 \\
  1 & \frac{N^\frac{3}{2}}{2} & 0 & -\frac{N^\frac{3}{2}}{2} & -1 \\
  1 & 0 & -2 & 0 & 1 \\
  1 & -\sqrt{N} & 0 & \sqrt{N} & -1 \\
  1 & -N & 2(N-1) & -N & 1 \\
    \end{array}
  \right)$ ,\\

$Q=\left(\begin{array}{ccccc}
 1 & N & \frac{(N-1)(N+2)}{2} & \frac{N^2}{2} & \frac{N}{2} \\
 1 & \sqrt{N} & 0 & -\sqrt{N} & -1 \\
 1 & 0 & -\frac{N}{2}-1 & 0 & \frac{N}{2} \\
 1 & -\sqrt{N} & 0 & \sqrt{N} & -1 \\
 1 & -N & \frac{(N-1)(N+2)}{2} & -\frac{N^2}{2} & \frac{N}{2} \\
 \end{array}
 \right)$ .\\

\bigskip
The intermediate association scheme $Z$ between $X$ and $Y$, and of size
$N^2/2$ in $\mathbb{R}^{N-1}$, where $Z$ is obtained from vectors
$Z'=\{x\in X\mid (x,u)=1/\sqrt{N}\}$ by orthogonal projection to
$\langle u\rangle ^\perp$ for any fixed $u\in X$ and rescaling to make
a spherical code, is a class 3 association scheme.
This is a spherical $3$-design and of degree $s=3$ with the following
parameters. It is interesting to note that they are $Q$-polynomial
association schemes
(and not $P$-polynomial association scheme for $N\geq 4$.)
It seems that this example was already recognized in the list of
W. Martin's homepage (of such association schemes) as those coming
from the linked symmetric designs.
These association schemes are also possible candidates of universally
optimal codes in the sense of Cohn-Kumar \cite{Cohn2}.
Here we describe the parameters of the association scheme $Z$.\\

$B_1=\left (\begin{array}{cccc}
 0 & 1 & 0 & 0  \\
\frac{(N-\sqrt{N})(N-2)}{4}  & \frac{(N-3\sqrt{N})(N-4)}{8} & \frac{(N-\sqrt{N})(N-4)}{8} & \frac{(N-2\sqrt{N})(N-2)}{8}\\
0 & \frac{(N+\sqrt{N})(N-4)}{8} & \frac{(N-\sqrt{N})(N-4)}{8} & \frac{N(N-2)}{8}   \\
0 & \frac{(\sqrt{N}-2)(\sqrt{N}+1)}{2} & \frac{N-\sqrt{N}}{2} & 0 \\
\end{array}
\right)$, \\

$B_2=\left (\begin{array}{cccc}
0 & 0 & 1 & 0 \\
0 & \frac{(N+\sqrt{N})(N-4)}{8} & \frac{(N-\sqrt{N})(N-4)}{8} & \frac{N(N-2)}{8} \\
\frac{(N+\sqrt{N})(N-2)}{4} & \frac{(N+\sqrt{N})(N-4)}{8} & \frac{(N+3\sqrt{N})(N-4)}{8} & \frac{(N+2\sqrt{N})(N-2)}{8} \\
0 & \frac{N+\sqrt{N}}{2} & \frac{(\sqrt{N}-1)(\sqrt{N}+2)}{2} & 0 \\
\end{array}
\right)$ ,\\

$B_3=\left (\begin{array}{cccc}
0 & 0 & 0 & 1 \\
0 & \frac{(\sqrt{N}-2)(\sqrt{N}+1)}{2} & \frac{N-\sqrt{N}}{2} & 0 \\
0 & \frac{N+\sqrt{N}}{2} & \frac{(\sqrt{N}-1)(\sqrt{N}+2)}{2} & 0 \\
N-1 & 0 & 0 & N-2 \\
\end{array}
\right)$ ;\\

$B_1^*=\left (\begin{array}{cccc}
 0 & 1 & 0 & 0  \\
 N-1 & 0 & 2 & 0 \\
 0 & N-2 & N-4 & N-1  \\
 0 & 0 & 1 & 0 \\
\end{array}
\right)$ ,\\

$B_2^*=\left (\begin{array}{cccc}
 0 & 0 & 1 & 0  \\
 0 & N-2 & N-4 & N-1 \\
 \frac{(N-2)(N-1)}{2} & \frac{(N-4)(N-2)}{2} & \frac{N^2-6N+12}{2} & \frac{(N-4)(N-1)}{2}  \\
 0 & \frac{N}{2}-1 & \frac{N}{2}-2 & 0 \\
\end{array}
\right)$ ,\\

$B_3^*=\left (\begin{array}{cccc}
 0 & 0 & 0 & 1  \\
 0 & 0 & 1 & 0 \\
 0 & \frac{N}{2}-1 & \frac{N}{2}-2 & 0  \\
 \frac{N}{2}-1 & 0 & 0 & \frac{N}{2}-2 \\
\end{array}
\right)$. \\

\medskip
\noindent
The first and the second eigenmatrices are the same as in Proposition \ref{short-dels}.

\medskip

\noindent
{\bf Remarks.}
(1) We have association schemes $X$, $Z$, $Y$ of sizes $N^2+2N$,
$N^2/2$, $N^2/4$ (respectively) in $\mathbb{R}^N$, $\mathbb{R}^{N-1}$,
$\mathbb{R}^{N-2}$ (respectively), where $N$ must be an even power of $2$.
Currently all of them are possible candidates of universally
optimal codes in the sense of Cohn-Kumar, at least for $N\geq 16$.
It is shown in Cohn et. al. \cite{Cohn3} that $X$ in $R^4$
(i.e., for $N=4$) is not universally optimal.
(It is an open question whether it is optimal or not.)
$Z$ for $\mathbb{R}^3$ (i.e., for $N=4$) is not universally optimal nor optimal.
On the other hand, $Y$ for $\mathbb{R}^2$ (i.e., for $N=4$)
is universally optimal. Although it is a wild guess without firm ground,
we think $Y$ may be most likely to be universally optimal among $X$, $Z$, and $Y$. \\
(2) It is an interesting open question whether these association
schemes are uniquely determined by the parameters.
The uniqueness of $Y$ for $\mathbb{R}^{14}$ (i.e., for $N=16$) was obtained in \cite{Ban}.
On the other hand, the uniqueness of $X$ for $\mathbb{R}^{16}$ (i.e., for $N=16$)
was proved by Akio Nakamura \cite{Nak} in his masters degree
thesis of Kyushu University in 1997 (it follows also from the uniqueness of
the Nordstrom-Robinson code).
We note that the uniqueness of $Z$ for $\mathbb{R}^{15}$ (i.e., for $N=16$)
is also obtained. The claim is essentially obtained in \cite{Mat}, it can
be proved also by method of \cite{Nak}.
So, it would be interesting what will happen in particular for $N=64$ for
$X$, $Z$, and $Y$. The result of Kantor \cite{Kan2} implies
that if $N=2^{m+1}$ with odd $m$, and if $m$ is not a prime,
then there are non-isomorphic line systems, and so there are non-isomorphic
association schemes with the same parameters, i.e., the uniqueness is break down. \\
(3) Quite recently it was shown \cite{Boy2} that the problem of constructing of
$s$ pairwise mutually unbiased bases in $\mathbb{K}^n$
($\mathbb{K}=\mathbb{R}$ or $\mathbb{K}=\mathbb{C}$) is equivalent to
the problem of constructing of $s$ Cartan subalgebras of $sl_n(\mathbb{K})$
that are pairwise orthogonal with respect to Killing form and are closed
under the adjoint operation.
In particular, a complete collection of mutually unbiased bases in
$\mathbb{C}^n$ is equivalent to an orthogonal decomposition of Lie algebra
$sl_n(\mathbb{C})$, closed under the adjoint operation.
So there is a link to the well-developed theory \cite{Kos} of orthogonal
decompositions of Lie algebras.

\section{Extremal line-sets and Barnes-Wall lattices}
\label{ext}

In this section we discuss connections between mutually unbiased bases,
extremal line-sets in $\mathbb{R}^N$ with prescribed angles
and minimum vectors of Barnes-Wall lattices.

Fix any positive integer $N>1$. Let $M$ be a set of unit vectors in
$\mathbb{R}^N$, such that $|(a,b)|\in \{0,1/\sqrt{N}\}$ for all $a\neq b$
in $M$ (so, in particular, $M \cap (-M) =\emptyset$).
Then $|M|\leq N(N+2)/2$, and if $|M|$ reaches this upper bound, then $M$ is
a set of $N/2 +1$ mutually unbiased bases \cite[Proposition 3.12]{Cal}.
Constructions of such extremal line-sets are known only for $N=2^{m+1}$,
$m$ odd \cite{Cal}.

Our final observation is that known constructions of extremal line-sets
(or, equivalently, maximal sets of mutually unbiased bases) are connected to
the minimum vectors of Barnes-Wall lattices.
We show that vectors of known maximal real MUB after suitable rescaling will
become minimal vectors of a Barnes-Wall lattice.
Therefore, vectors of association schemes $X$, $Y$, $Z$ can be obtained
from a set of minimal vectors of the Barnes-Wall lattices.

First we recall the construction from \cite{Cal}. Label the standard basis of
$\mathbb{R}^N$ as $e_v$, with $v\in V=\mathbb{Z}_2^{m+1}$.
For $b\in V$, define the permutation matrix $X(b)$ and diagonal matrix
$Y(b)$ as follows:
$$X(b):e_v\mapsto e_{v+b} \ \ {\rm and} \ \ Y(b):= {\rm diag} [(-1)^{b\cdot v}]. $$
The groups $X(V):=\{X(b) \mid b\in V\}$ and $Y(V):=\{Y(b) \mid b\in V\}$
are contained in $O(\mathbb{R}^N)$ and are isomorphic to the additive
group $V$. Let $E:=\langle X(V),Y(V)\rangle$.
Then $E=2_+^{1+2(m+1)}$ is an extraspecial 2-group of
order $2^{1+2(m+1)}$ and $\overline{E}=E/Z(E)$ is elementary abelian
group of order $2^{2(m+1)}$.  We identify the center $Z(E)$ of $E$ with
$\mathbb{Z}_2$ and consider the map $Q: \overline{E}\rightarrow \mathbb{Z}_2$
defined by $Q(\overline{e})=e^2$ for any $\overline{e}\in \overline{E}$
and any preimage $e$ of $\overline{e}$ in $E$. Then $Q$ is a non-singular
quadratic form on $\overline{E}$.
So $\overline{E}$ is an $\Omega^+(2m+2,2)$-space. The action of $E$
on $\mathbb{R}^N$ can be extended to the action of the group
$2_+^{1+2(m+1)}\Omega^+(2m+2,2)$.

The space $\overline{E}$ contains $(2^{m+1}-1)(2^m+1)$ singular points.
An orthogonal spread of $\overline{E}$ is a family $\Sigma$ of
$2^m+1$ totally singular $(m+1)$-spaces such that every singular point of
$\overline{E}$ belongs to exactly one member of $\Sigma$.

Let $A$ be a subgroup of $E$ such that its image $\overline{A}$ in
$\overline{E}$ is totally singular $(m+1)$-space of $\overline{E}$.
Then the set $\mathcal{F}(A)$ of $A$-irreducible subspaces of $\mathbb{R}^N$
is an orthogonal frame: a set of $2^{m+1}$ pairwise orthogonal lines
through the origin. For an orthogonal spread $\Sigma$ of the
$\Omega^+(2m+2,2)$-space $\overline{E}$ we let
$$ \mathcal{F}(\Sigma):=\bigcup_{\overline{A}\in \Sigma} \mathcal{F}(A). $$
Then $\mathcal{F}(\Sigma)$ consists of $2^{m+1}(2^m+1)$ lines of $\mathbb{R}^N$
such that, if $u_1$ and $u_2$ are unit vectors in different members of
$\mathcal{F}(\Sigma)$, then $|(u_1,u_2)|=0$ or $2^{-(m+1)/2}$.
Therefore, $\mathcal{F}(A)$ determines orthonormal basis and
$\mathcal{F}(\Sigma)$ determines a set of $N/2+1$ mutually
unbiased bases.
These line-sets $\mathcal{F}(\Sigma)$ are extremal in the sense that
$|\mathcal{F}(\Sigma)|$ meets an upper bound obtained in \cite{Del2}
for line-sets in $\mathbb{R}^N$ with prescribed angles.

The binary Kerdock code $\mathcal{K}(\Sigma)$ can be recovered
\cite{Cal} from $\mathcal{F}(\Sigma)$:
$$\mathcal{K}(\Sigma)= \{(c_v)_v\in \mathbb{Z}_2^N \mid
\langle((-1)^{c_v})_v\rangle\in \mathcal{F}(\Sigma)\}.$$

\medskip

Take unit vectors from line-set $\mathcal{F}(\Sigma)$, rescale them to
vectors of norm $\sqrt{N}$, then these vectors will be minimum
vectors of a Barnes-Wall lattice.
Indeed, we note that for odd $m$ the minimum norm of Barnes-Wall
lattice is $\sqrt{N}$, the automorphism group is
$G=2^{1+2(m+1)}\Omega^+(2m+2,2)$, $G$ acts transitively \cite{Gri}
on the set of minimum vectors, the vector $c=N^{-1/4}\sum_{v\in V}e_v$
is a minimum vector, and any even lattice of rank $N$ invariant under the group
$G$ is similar to a Barnes-Wall lattice \cite{Gri}.
So all the minimum vectors of Barnes-Wall lattice are obtained
from $c$ by action of the group $G$.
On the other hand, in notations of \cite{Cal} we have $c=N^{1/4}e_b^*$
for $b=0$, and any unit vector of $\mathcal{F}(\Sigma)$ is obtained
from $e_b^*$ by action of some element of $G$.



\begin{thebibliography}{99}

\bibitem{Bal}
B. Ballinger, G. Blekherman, H. Cohn, N. Giansiracusa, E. Kelly and A. Sch\"urmann,
Experimental study of energy-minimizing point configurations on spheres,
arXiv:math/0611451v2 [math.MG].

\bibitem{Ban}
E. Bannai, E. Bannai and H. Bannai,
Uniqueness of certain association schemes,
to appear in {\em European Journal of Combinatorics}.

\bibitem{Ban2}
E. Bannai and E. Bannai,
On antipodal spherical $t$-designs of degree $s$ with $t\geq 2s-3$,
preprint.

\bibitem{Boy}
P. O. Boykin, M. Sitharam, M. Tarifi, P. Wocjan,
Real mutually unbiased bases, arXiv:quant-ph/0502024v2.

\bibitem{Boy2}
P. O. Boykin, M. Sitharam, Pham Huu Tiep and P. Wocjan,
Mutually unbiased bases and orthogonal decompositions of Lie algebras,
arXiv:quant-ph/0506089v1.

\bibitem{Caen} D. de Caen and E. R. van Dam,
Association schemes related to Kasami codes and Kerdock sets,
Designs and codes - a memorial tribute to Ed Assmus.
{\em Des. Codes Cryptogr.} {\bf 18} (1999), no. 1--3, 89--102.

\bibitem{Cal}
A. R. Calderbank, P. J. Cameron, W. M. Kantor and J. J. Seidel,
$\mathbb{Z}_4$-Kerdock codes, orthogonal spreads, and extremal Euclidean line-sets,
{\em Proc. London Math. Soc.} (3) 75 (1997), 436--480.

\bibitem{Cohn}
H. Cohn,
Sphere packings, energy minimization, and linear programming bounds,
in The Proceedings of Second COE Workshop on Sphere Packings, (2005), 1--42.

\bibitem{Cohn2}
H. Cohn and A. Kumar,
Universally optimal distributions of points on spheres,
{\em J. Amer. Math. Soc.} 20, No. 1 (2007) 99--148.

\bibitem{Cohn3}
H. Cohn, J. H. Conway, N. D. Elkies and A. Kumar,
The $D_4$ root system is not universally optimal,
to appear in {\em Experimental Mathematics}, arXiv:math/0607447v2 [math.MG].

\bibitem{Cam}
P. Camion,
Codes and association schemes: Basic properties of
association schemes relevant to coding,
In "Handbook of Coding Theory", (V. S. Pless and W. C. Huffman, Eds.)
Volume 2, Chapter 18. Elsevier, Amsterdam, 1998.

\bibitem{Del}
P. Delsarte,
An algebraic approach to the association schemes of coding theory,
{\em Phillips Res. Repts. Suppl.} 10 (1973).

\bibitem{Del2}
P. Delsarte, J. M. Goethals and J. J. Seidel,
Bounds for systems of lines and Jacobi polynomials,
{\em Phillips Res. Repts.} 30 (1975) 91--105.

\bibitem{Del3}
P. Delsarte, J. M. Goethals and J. J. Seidel,
Spherical codes and designs,
{\em Geom. Dedicata} 6 (1977) 363--388.


\bibitem{Gri}
R. L. Griess Jr.,
Pieces of $2^d$: existence and uniqueness for Barnes-Wall lattices
and Ypsilanti lattices.
{\em Advances in Mathematics}, 196 (2005) 147--192.

\bibitem{Gri2}
R. L. Griess Jr.,
Few-cosine spherical codes and Barnes-Wall lattices,
arXiv: math/0605175v1 [math.CO].

\bibitem{Ham}
R. Hammons, P. V. Kumar, A. R. Calderbank, N. J. A. Sloane and P. Sol\'{e},
The $Z_4$-linearity of Kerdock, Preparata, Goethals, and related codes,
{\em IEEE Trans. Inform. Theory} 40, No. 2 (1994), 301--319.

\bibitem{Kan}
W. M. Kantor,
Spreads, translation planes and Kerdock sets, I, II,
{\em SIAM J. Alg. Discr. Math.} 3 (1982) 151--165, 308--318.

\bibitem{Kan2}
W. M. Kantor,
An exponential number of generalized Kerdock codes,
{\em Inform. Control} 53 (1982) 74--80.

\bibitem{Kos}
A. I. Kostrikin and Pham Huu Tiep,
Orthogonal decompositions and integral lattices,
Walter de Gruyter, Berlin (1994).

\bibitem{Mat}
R. Mathon,
The systems of linked $2-(16,\,6,\,2)$ designs,
{\em Ars Combin}. 11 (1981) 131--148.


\bibitem{Nak}
A. Nakamura,
On extremal line set in Euclidean space with prescribed angles,
Masters degree thesis, Kyushu University, 1996.



\end{thebibliography}
\end{document}